\documentclass[review]{elsarticle}
\usepackage{lineno}
\usepackage[pdftex]{hyperref}
\modulolinenumbers[5]

\journal{Journal of Computational and Applied Mathematics}

\usepackage{amsmath}
\usepackage{amssymb}
\usepackage{graphicx}

\newcommand{\bx}{{\bf x}}
\newcommand{\I}{\text{i}}
\newtheorem{theorem}{Theorem}
\begin{document}
\begin{frontmatter}
\title{Using B-spline frames to represent solutions of acoustics scattering problems}

\author[adressARI]{Wolfgang Kreuzer}
\ead{wolfgang.kreuzer@oeaw.ac.at}
\address[adressARI]{Acoustics Research Institute, Austrian Academy of Sciences, \\Wohllebengasse 12-14, Vienna, 1040, Vienna}

\begin{abstract}
Although frames, which are a generalization of bases, are important tools used in signal processing, their potential in other fields of engineering and applied mathematics (e.g. acoustics) has not been fully explored yet. Gabor frames, that are specific type of frames, are very well adapted to oscillating functions, and therefore have a great potential to efficiently represent functions in connection with the Helmholtz operator. In this paper representations of the solution of a scattering problem  in 2D using Gabor frames based on B-splines as building blocks are investigated. Practical issues concerning the implementation of frames like the restriction of the frame elements to a finite interval and  methods to determine the unknown coefficients for the representation, i.e. by using dual frames or by solving a least squares problem, are discussed. Numerical experiments are made comparing the different ways to determine the unknown expansion coefficients and  the representations are compared in terms of their efficiency, i.e. the number of used frame coefficients versus the accuracy of the approximation. In all cases the coefficients calculated with a slightly modified orthogonal matching pursuit algorithm provide the best accuracy versus sparsity ratio.
\end{abstract}


\end{frontmatter}


\section{Introduction}
Numerical methods like the boundary element method (BEM) or finite elements (FEM) play an important role for calculating the scattering of acoustic waves from complex structures, thus for solving the Helmholtz equation (cf. \cite{Gauletal03,Fischeretal04,SauSch11,Chenetal08}). To that end, the scatterer (or its surface) is discretized with simple geometric elements and based on this discretization the unknown solution $u(\bx)$ (i.e. the acoustic pressure, the velocity potential or the particle velocity) is represented/approximated using simple (basis) functions: $u(\bx) = \sum_{n=0}^N u_n \phi_n(\bx)$. In general, for Helmholtz problems with uniform grids the gridsize and thus the number of necessary ansatz functions $N$ is dependent on the wavenumber $k=\frac{2\pi f}{c}$ where $f$ is the frequency and $c$ the speed of sound (cf. \cite{Marburg02,Langeretal17}). Therefore, for high frequencies the computational effort becomes very large, and calculations of wave scattering problems can last between hours to several days on modern computers.

One idea to circumvent the frequency dependence of the mesh (e.g. for BEM six to eight elements per wavelength is recommended, for FEM similar rules exist) is to include oscillating components in the ansatz functions (e.g. see \cite{BruGeu07,Chandler-Wildeetal12}). As the construction of a basis with special properties can be cumbersome sometimes, the potential of ansatz functions based on \emph{Gabor frames} \cite{FeiStr98,Groechenig01,Christensen08,Christensen03}, i.e. functions of the type $\phi_{m,n}(x) = g(x - na)e^{2\pi\I mb x}, m,n \in \mathbb{Z}$, is investigated in this work. Frames are generalizations of bi-orthogonal bases and, contrary to bases, lead to possibly redundant 
representations. Because  of  the relaxed requirements it is easier to construct frames with special a-priori properties (e.g. a sparse representation) compared to finding appropriate bases (for a hands-on survey on frames and redundancy cf. \cite{KovChe07a,KovChe07b}).  A motivation for the advantage of using a (redundant) frame instead of a  basis can be given by the following example: In a ``rich'' dictionary with a lot of entries, it is much easier to find the correct pieces to have a short, i.e. sparse, representation of a given sentence. 

Similar to Riesz bases, frames allow to represent arbitrary elements of a given Hilbert space by series expansions with respect to the frame elements. Frames have already found some application in signal processing \cite{Boelcskeietal98,Balazsetal13} and psychoacoustics \cite{Balazs10,Balazsetal17}. 

In this paper the idea of using frames for approximating solutions of scattering problems is presented. There is already some literature in connection with frames and operator equations, e.g. \cite{Stevenson03, Dahlkeetal07, Dahlkeetal07a, Balazs08,Harbrechtetal08}, but most papers are rather on a theoretical and conceptional level. This manuscript aims at a more practical and applied viewpoint of this topic. Some aspects of implementing and applying frames are discussed, and numerical experiments for two different wavenumbers are performed. In this paper, the problems discussed are restricted to 2D scattering problems, however, the extension to higher dimensions is straight forward by using tensor products of one-dimensional frames. It is clear, that for solving 2D scattering problems already very efficient methods exist (e.g. spectral methods based on Nystr\"om methods \cite{ColKre13}), however, in this paper the focus lies on investigating the potential that frames offer outside the area of signal processing especially with respect to providing efficient representations of oscillating solutions of scattering problems. In that respect, this manuscript should be seen as one of the first steps towards efficient methods for solving 3D scattering problems with relatively large wavenumbers. The main aim of this manuscript is to introduce the concept of frames to the field of applied sciences away from signal processing and to address the question: \emph{``Can Gabor frames be used to provide an efficient representation for solutions of scattering problems?''}.

The paper is structured as follows. In Section \ref{Sec:Bsplines} a brief overview of B-spline functions, that will play the role as generating window for the Gabor frames, is given. In Section \ref{Sec:Frames} the definition of (Gabor)frames and their duals is given and some of their properties are discussed, especially in connection with Gabor frames based on B-splines as generating window functions. Section \ref{Sec:PracticalAspects} deals with practical aspects for implementing frames, e.g. restriction to a finite interval, strategies for sampling and ways to calculate the expansion coefficients for arbitrary functions in $L^2(\mathbb{R})$. In Section \ref{Sec:Experiments} numerical experiments dealing with the representation of the solution of a scattering problem in 2D using B-spline Gabor frames will be performed. In this section the unknown expansion coefficients of the frame representation are calculated either by using the  product of three known dual frames with the target functions or by solving a least squares problem using the orthogonal matching pursuit algorithm (OMP). The OMP algorithm was slightly modified so that in one iteration more than one atom (= frame element) can be chosen and the effect of the number of atoms chosen per iteration on the frame coefficients are investigated.  The numerical experiments are performed  for two different frequencies and the representations are compared with respect to their efficiency in terms of accuracy of the representation versus the number of coefficients used. As target functions for the approximations the sound field on a sound hard circle caused by a point source outside the circle and the Hankel function of order 0, which is the free space Green's function for the 2D Helmholtz equation,  are chosen. 

The octave scripts used in the numerical experiments can be found at \url{https://www.kfs.oeaw.ac.at/research/projects/biotop/Bsplineframes.tgz}, for a description of the different scripts refer to the file README.txt in the tar-ball.

\section{B-splines}\label{Sec:Bsplines}
\begin{figure}[!htb]
  \begin{center}
  \includegraphics[width=0.5\textwidth]{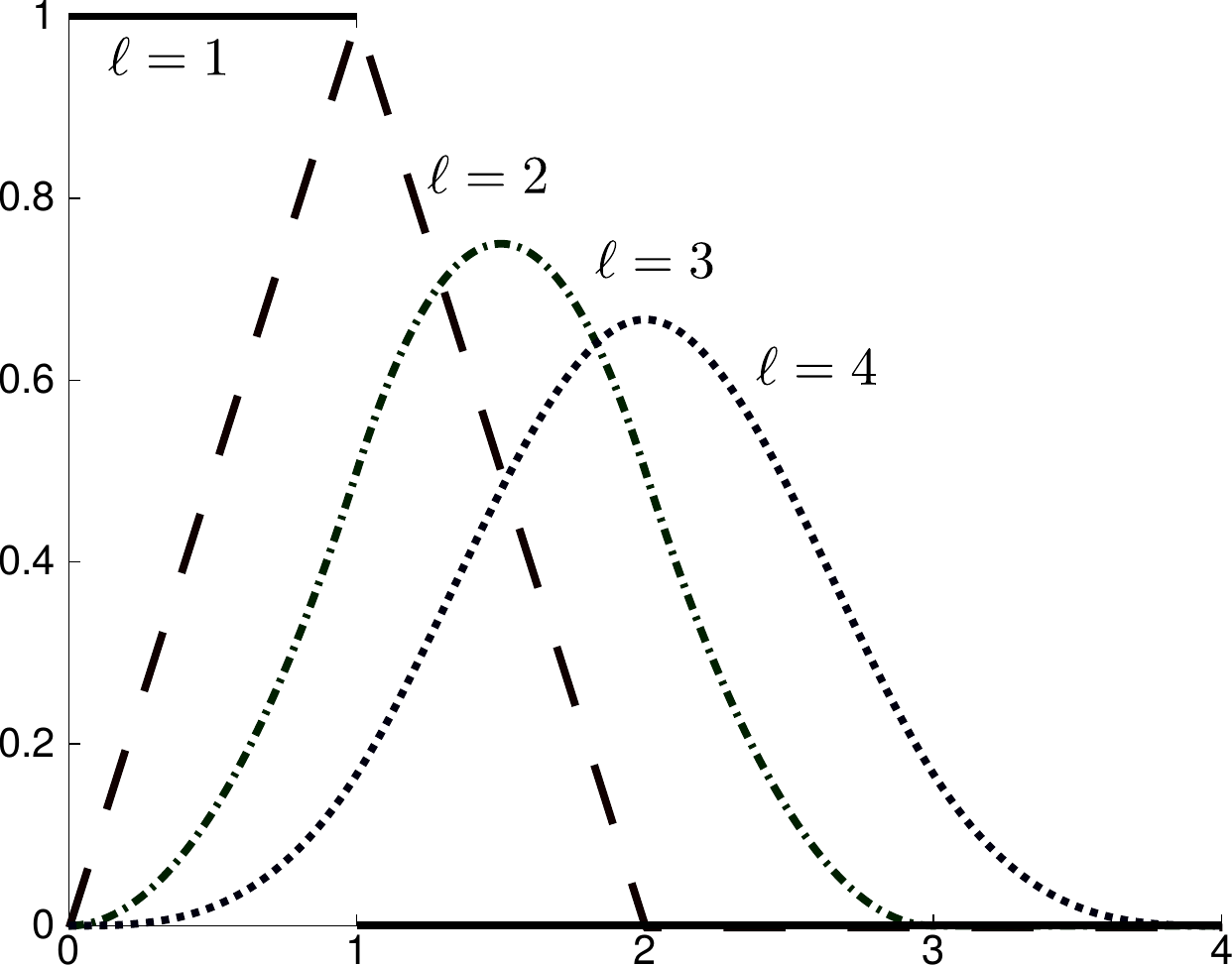}
  \caption{B-splines of order one to four}\label{Fig:Bsplines}
  \end{center}
\end{figure}

B-spline functions play an important role in numerical and applied mathematics (e.g. as ansatz functions in FEM and BEM or in connection with NURBS curves and isogeometric analysis \cite{Hughesetal10,Doelzetal18}) and especially in computer graphics (e.g. \cite{deBoor78,Bartelsetal86,PieTil97}). They can be easily constructed via
\begin{align}
  N_1(x) = \left\{
  \begin{array}{c}
    1,\quad x\in[0,1],\\
    0 \quad \text{otherwise},
  \end{array}
  \right.\\
  N_{\ell+1}(x) = (N_\ell* N_1)(x) = \int_0^1 N_\ell(x-t)dt,
\end{align}
where '*' denotes the convolution, for examples refer to Fig.~\ref{Fig:Bsplines} for the B-splines $N_\ell$ for orders $\ell = 1,\dots,4$.
Besides their finite support B-splines offer additional interesting properties, e.g.
\begin{itemize}
\item $\int\limits_{-\infty}^{\infty} N_\ell(x)dx = 1$,
\item partition of unity, i.e.  $\sum\limits_{k\in\mathbb{Z}} N_\ell(x-k) = 1$,
\item $N_\ell(x) \in C^{\ell-2}(\mathbb{R})$, i.e. $N_\ell$ is $\ell-2$ times continuously differentiable,
\item $N_\ell$ restricted to $[k,k+1]$ is a polynomial of degree at most $\ell-1$ for all $k\in \mathbb{Z}$,
\item $N_\ell(.+k),\; k = 0, \dots \ell-1$ are linearly independent on $[0,1]$,
\item $\int\limits_{-\infty}^\infty N_\ell(x)f(x) dx = \int\limits_{[0,1]^\ell} f(x_1 + \cdots + x_\ell) d_{x_1}\cdots d_{x_\ell}$.
\end{itemize}
The linear independence, the partition of unity property, and the compact support make (low order) B-splines good candidates for FEM and BEM ansatz functions, and in our case also as generating window functions for Gabor frames.
\section{Frames}\label{Sec:Frames}
 A countable family of functions $\{g_i\}_{i\in \mathcal{I}}$ in a Hilbert space $\mathcal{H}$ is called \emph{frame} for $\mathcal{H}$ if there exist constants $A,B > 0$ such that
\begin{equation}\label{Equ:FrameCondition}
  A||f||^2 \le \sum_{i\in \mathcal{I}} | \langle f,g_i \rangle |^2 \le B||f||^2,\quad \forall f \in \mathcal{H},
\end{equation}
where $\mathcal{I} \subseteq \mathbb{N}$ is some countable index set. $\langle f,g\rangle$ denotes the scalar product in the Hilbert space, e.g. for $f,g \in L^2(\mathbb{R})$ the product is given by $\langle f,g \rangle = \int_{-\infty}^\infty f(x)g^*(x)dx$, where $g^*(x)$ denotes the complex conjugate of the function $g(x)$.

The bounds $A$ and $B$ are called \emph{frame bounds}, if $A = B$ the frame is called  \emph{tight}. In the Hilbert space of square integrable functions $L^2\left(\mathbb{R}\right)$ Eq.~(\ref{Equ:FrameCondition}) can be seen as the link between the norm/energy of a function and the norm/energy of its representation using frame atoms $g_i(x)$. Eq.~(\ref{Equ:FrameCondition}) ensures that every element in the vector space can be reconstructed in a stable way using the frame atoms, if the frame bounds are close to each other the reconstruction is faster and behaves numerically better.

For every frame $\{g_i\}_{i\in \mathcal{I}}$ there exists at least one \emph{dual} frame $\{\tilde{g}_i\}_{i\in \mathcal{I}}$ such that
\begin{equation}\label{Equ:FrameRep}
  f = \sum_{i\in \mathcal{I}}\langle f,\tilde{g}_i \rangle g_i = \sum_{i\in \mathcal{I}}\langle f,g_i \rangle \tilde{g}_i,
\end{equation}
where the sums in the above equation converge unconditionally. Thus, every element in the Hilbert space can be represented by a weighted (possible infinite) sum of frame atoms and the coefficients of this representation can be calculated by the inner product of the target function with the dual frame. A frame is similar to a basis, but as frames can have multiple different dual frames the expansion is \emph{not} unique. Amongst all possible dual frames the \emph{canonical dual} frame plays a special role as it can be constructed by inverting the so called \emph{frame operator} $f\rightarrow \sum_{i\in\mathcal{I}} \langle f,g_i\rangle g_i$, which in practice means by taking the pseudo inverse of the matrix containing the sampled frame elements (see also Section \ref{Sec:PracticalAspects}). For a more detailed introduction to frames refer for example to \cite{Christensen01,Christensen08}.
\subsection{Gabor frames}\label{Sec:Gabor}
In the Hilbert space of square integrable functions $L^2(\mathbb{R})$  Gabor systems $\mathcal{G}(g,a,b)$ are collections of functions  that are constructed by translating and modulating a given window function $g(x)$:
\begin{equation}
  g_{mn}(x) = E_{mb}T_{na} g(x) = g(x-na)e^{2\pi\I mbx}, 
\end{equation}
where $a,b \in \mathbb{R}^+$ and $m,n\in \mathbb{Z}$. $E_{mb} f(x) = f(x) e^{2\pi\I mbx}$ and $T_{na} f(x) = f(x - na)$ define the modulation and the translation operators, respectively. Under certain conditions on the parameters $a$ and $b$ a Gabor system forms a frame. For Gabor frames it is known that some dual frames (especially the canonical dual) also have a Gabor structure, thus the duals can be constructed by translation and modulation of a \emph{dual window} function (c.f. \cite{Perraudinetal13,Strohmer98}). 

A continuous version of a Gabor frame is given by the short time Fourier transform
\begin{equation}
  \text{STFT}\{s(t)\}(\tau,\omega) = \int_{-\infty}^\infty s(t) g(t-\tau)e^{-\I \omega t} dt,
\end{equation}
where $g(t)$ is a specific window function (e.g. Hanning window). The STFT is essential for the time-frequency representation of signals $s(t)$, e.g. for generating spectrograms and Gabor frames lead to a sampled version of the STFT. 

\subsection{Gabor frames based on B-splines}\label{Sec:Bsplineframes}
Because of their compact support and especially because of their partition of unity property B-splines are easy to use window functions for generating Gabor frames and some of their duals. Based on the theorems given by Christensen \cite{Christensen03,Christensen08,Christensen10} conditions for the frame parameters $a,b$ can be found to derive Gabor frames based on B-splines \emph{and} to construct some duals frames:
\begin{theorem}\label{Th:Gen}
  For $\ell \in \mathbb{N}$, the B-spline $N_\ell(x)$ generates Gabor frames for all $(a,b) \in (0,\ell]\times(0,1/\ell]$.
\end{theorem}
\begin{theorem}\label{Th:Dual1}
  For any $\ell \in \mathbb{N}$, and $b\in (0,\frac{1}{2\ell-1}]$, the functions $N_\ell$ and
  \begin{equation}\label{Equ:Dual1}
    h_\ell = b N_\ell(x) + 2b\sum_{k=1}^\ell N_\ell(x+k)
  \end{equation}
  generate dual frames $\{E_{mb}T_nN_\ell\}_{m,n\in \mathbb{Z}}$ and $\{E_{mb}T_nh_\ell\}_{m,n\in \mathbb{Z}}.$
\end{theorem}
\begin{theorem}\label{Th:Dual2}
  For any $\ell \in \mathbb{N}$ and $b\in (0,\frac{1}{2\ell - 1}]$, the functions $N_\ell$ and
  \begin{equation}\label{Equ:Dual2}
    h_\ell(x) = b\sum_{k = -\ell + 1}^{\ell -1}N_\ell(x+k)
  \end{equation}
  generate dual frames $\{E_{mb}T_nN_\ell\}_{m,n\in \mathbb{Z}}$ and $\{E_{mb}T_nh_\ell\}_{m,n\in \mathbb{Z}}.$
\end{theorem}

\begin{figure}[!thb]
  \begin{center}
    \includegraphics[width=0.5\textwidth]{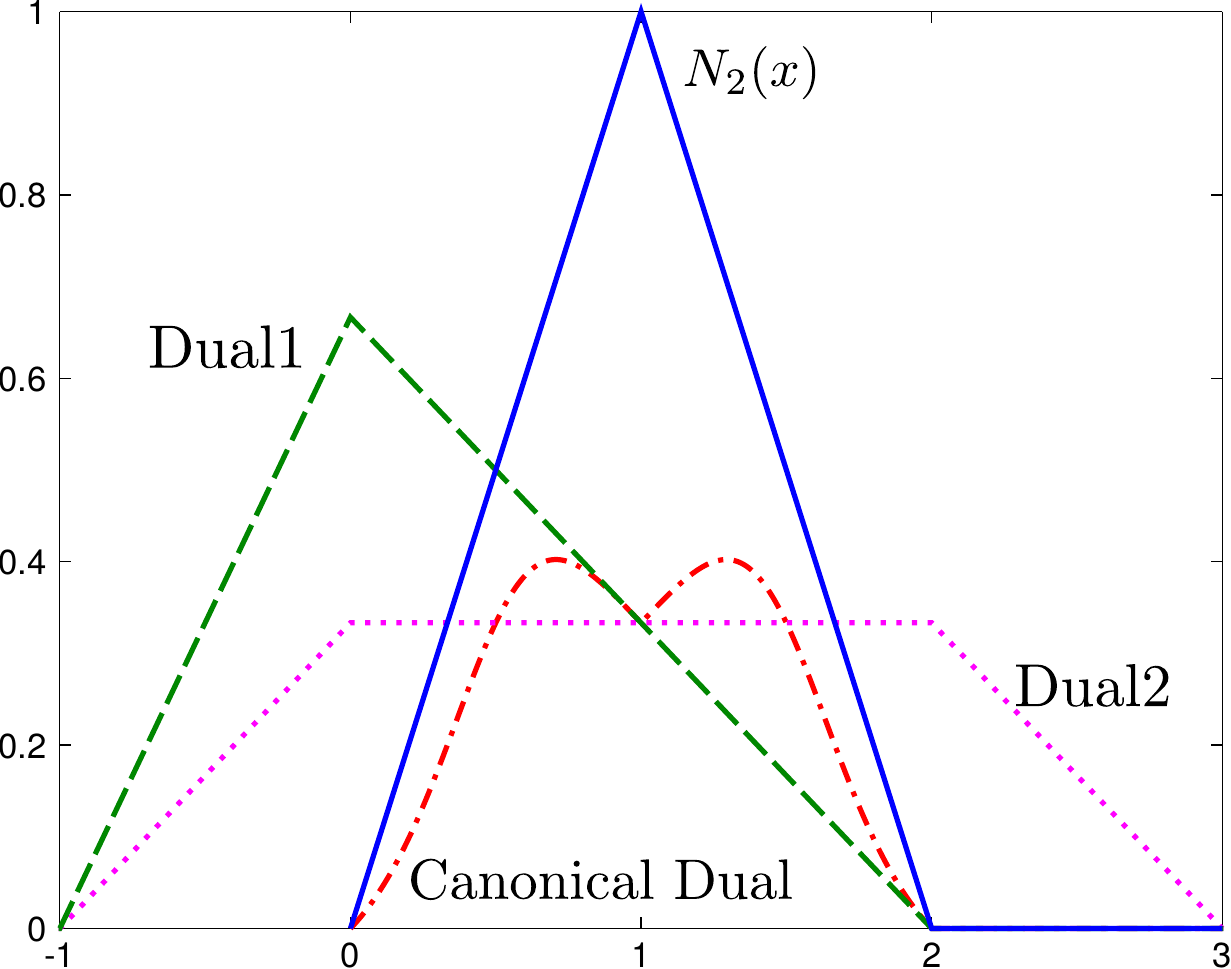}
    \caption{B-spline of order 2 (solid line) and dual windows created by applying Eqs.~(\ref{Equ:Dual1}) and (\ref{Equ:Dual2}) for $a = 1$ and $b = \frac13$. The dual windows are depicted using the dashed and dotted lines,  respectively. Additionally the canonical dual window (dash-dotted line) is depicted.}\label{Fig:Duals}
  \end{center}
\end{figure}

For the proofs refer to Corollaries 9.1.9 and 9.4.2 in \cite{Christensen08} and Corollary 2.46 in \cite{Christensen10}. In Fig.~\ref{Fig:Duals} the window function $N_2(x)$ and the two dual windows Dual1 and Dual2 (dashed and dotted, respectively) constructed  using Eqs.~(\ref{Equ:Dual1}) and (\ref{Equ:Dual2}) for $a = 1$ and $b = \frac13$ are shown. In Fig.~\ref{Fig:Duals} the (numerically determined) canonical dual window is shown with the dash-dotted line.  
\section{Practical Aspects}\label{Sec:PracticalAspects}
For practical application and implementation two points besides the construction of the frame have to be considered: a) The restriction to a finite interval, and b) sampling and the efficient calculation of the expansion coefficients.
\subsection{Restriction to Finite Intervals}\label{Sec:RestrictionInterval}
For a practical implementation it will be necessary to restrict the area of interest and thus the target function $f(x)$ to a finite interval, e.g. $[0,L]$. This can be done either by assuming $f(x)$ to be periodic with period $L$ or by ``ignoring'' all values of $f(x)$ outside $[0,L]$.  In signal processing periodization of the signal and the generating window is a standard option \cite{Sondergaard07}, however, for solutions of e.g. BEM or FEM problems, that are defined on finite patches and that are, in general, not periodic, this is sub optimal. In these cases periodization may lead to unwanted discontinuities at 0 and $L$ which can result in a high number of ``necessary'' frame elements around the interval boundaries. ``Necessary'' in this context means that 
in the representation many coefficients for frames located around the boundaries will have large absolute values, thus prohibiting an efficient expansion. 
However, setting all entries of $f(x)$ outside the interval of interest to zero also results in discontinuities at the interval boundaries,  which per se are not a problem for the frame used in the expansion, but for the multiplication with the dual frame. One way of dealing with this problem is to expand the interval of interest slightly so that it fully covers the support of the dual frame elements necessary for calculating the frame coefficients. An example: Suppose the function $f(x)$ should be represented in the interval [0,3] by using frames based on the B-spline of order 2 (the solid line in Fig.~\ref{Fig:Duals}). To fully cover the interval 4 shifted windows $N_2(x - n)$ with $n \in \{-1,0,1,2\}$ are needed. To calculate the frame coefficients using the dual frames the same amount of shifts of the dual window, which may have a bigger support (e.g. the dual window constructed using Theorem \ref{Th:Dual1} has support $[-1,3]$),  is necessary. Thus, for the correct calculation of the frame coefficients $f(x)$ is needed in the larger interval $[-2,5]$. If the function $f(x)$ is not known there, $f(x)$ should be expanded smoothly to that interval. For the reconstruction only parts of the frame elements inside $[0,3]$ will be used.
\subsection{Sampling}
In Sections \ref{Sec:Gabor} and \ref{Sec:Bsplineframes} it was assumed that all the frame elements are functions in $L^2(\mathbb{R})$. For numerical computations it is necessary to do some discretization at one point. This can be done in two ways. One way is to work in the function space and discretize the inner product Eq.~(\ref{Equ:FrameRep}) 
\begin{equation}\label{Equ:Quadrature}
  \langle f, \tilde{g}_i \rangle = \int_{-\infty}^\infty f(x) \tilde{g}^*_i(x) dx \approx \sum_{j = 1}^n \omega_j f(x_j)\tilde{g}_i^*(x_j), 
\end{equation}
where $\tilde{g}^*(x)$ denotes the complex conjugate of $\tilde{g}(x)$, and $x_j$ and $\omega_j$ are the nodes and weights of an appropriate quadrature method \cite{Evans94,Dominguezetal13,Dominguez14}, which is adapted to oscillating integrals. Since B-splines have finite support, all integrals in Eq.~(\ref{Equ:Quadrature}) are definite. 

Alternatively, one could already look at a discretized frame (cf. \cite{Sondergaard07} for an overview on sampled and periodized frames). In that case all inner products are simple vector products. Also the canonical dual frame can be calculated easily by taking the pseudo inverse of the matrix containing the sampled frame elements. For the discretized frame the number of frame elements will be finite because the number of necessary modulations is usually a function of the interval length, the frame parameter $b$, and the number of sampling points used. However, for discretized frames the frame parameters, the number of modulations, and the number of sampling points have to be chosen more carefully (cf. \cite{Sondergaard07}).  
\subsection{Calculation of the expansion coefficients}
If the dual frame is known (either as a function or in a discretized version) it is possible to calculate the coefficients $c_i = \langle f, \tilde{g}_i \rangle$ of the expansion $f(x) = \sum_{i = 1}^M c_i g_i(x)$ using the inner product of the target function $f(x)$ with the elements of the dual frames $\tilde{g}_i$. Alternatively, the unknown frame coefficients can be determined by finding the least squares solution
\begin{equation}\label{Equ:LeastSquares}
  \text{argmin}_{\bf c} ||f(x) - \sum_{i = 1}^N c_i g_i(x)||^2.
\end{equation}
In this case there is no need to know any dual frame, and the fact that a frame is used ensures that Eq.~(\ref{Equ:LeastSquares}) has a stable solution. However, finding a solution to Eq.~(\ref{Equ:LeastSquares}) may take some computation time, and because of the redundancy of a frame there will be several possible solutions. Solving Eq.~(\ref{Equ:LeastSquares}) using the pseudoinverse, which is equivalent of using the canonical dual frame, will, in general, not result in the most efficient solution in terms of sparsity. To arrive at an efficient representation methods from the field of compressed sensing should be used \cite{MalZha93,NeeTro09,EldKut12} e.g. the orthogonal matching pursuit (OMP) algorithm \cite{Patietal93,CaiWan11}.

\subsubsection{Orthogonal Matching Pursuit}\label{Sec:OMP}
The OMP algorithm is a simple greedy algorithm for finding a sparse solution of a least squares problem $\text{argmin}_{\bf c} ||{\bf Gc - f}||^2$. In case of a sampled discrete frame each column of the matrix ${\bf G}$ contains one frame element sampled at nodes $\bx_j, j = 1\dots N$, thus ${\bf G}_{ji} = g_i(\bx_j)$ and ${\bf f}_j = f(\bx_j)$, where $g_i(\bx_j)$ is the $i$-th frame element sampled at $\bx_j$. ${\bf G}^H$ is the Hermitian transpose of ${\bf G}$.

For the numerical experiments a slightly modified version of the (vector based) OMP MATLAB code from Stephen Becker \cite{Becker16} which partly follows the algorithm proposed in \cite{CaiWan11,EldKut12,NeeTro09} was used. The code is used in the ``straightforward'' slow mode, that does no orthogonalization of the chosen atoms as it was originally described for example in \cite{Patietal93}, but calculates the unknown coefficients by directly solving a least squares problem using the 'backslash' operator in MATLAB. In its most simple form the algorithm consists of the following 6 steps:
\begin{enumerate}
\item ${\bf r} = {\bf f}$.
\item Find the right candidate: $i_0 = \text{argmax} |{\bf G}^H {\bf r}|$.
\item Add the candidate to the list of used frame atoms: $\mathcal{I} = \mathcal{I} \cup \{i_0\}$.
\item Find ${\bf c}_{\mathcal{I}} = \text{argmin}_{\bf c} ||{\bf f} - {\bf G}_{\mathcal{I}} {\bf c}||^2$ where ${\bf G}_{\mathcal{I}}$ only contains the columns of ${\bf G}$ with indices in $\mathcal{I}$. 
\item Update the residual ${\bf r} = {\bf f} - {\bf G}_{\mathcal{I}} {\bf c}_{I}$.
\item If $\text{number}_\text{iterations} < \text{max}_\text{iterations}$ jump to step 2,
\end{enumerate}
where the number of maximum iterations has to be provided by the user.
As an alternative to the last step it is also possible to stop the iterations if the norm of the residuum is below a given  tolerance $tol$ provided by the user.

For non-sampled frames the OMP can be adapted to handle functions. In that case the multiplication in Step 2 has to be replaced by an appropriate quadrature method to solve $\langle f,g_i \rangle$ and the least squares problem has to be solved using the Gram matrix ${\bf G}_{ij} = \langle g_i,g_j \rangle, i,j \in \mathcal{I}$ (see \ref{Sec:OMPL2} for more details). The quadrature method used should be adapted to oscillating functions, for calculating the Gram matrix one should also use the fact that the B-splines are polynomials and that calculations can be done analytically by partial integration.

In the literature one can find several modifications of the OMP algorithm (cf. \cite{Blumensathetal12} for a small overview of some variants), in the numerical experiments below also the performance of a simple block-OMP algorithm will be investigated, where the frame elements corresponding to the $n$ largest entries of  $|{\bf G}^H {\bf r}|$ are chosen as candidates, with $n$ being the \emph{blocksize}. Judging from the experiments below the block-OMP has the advantage that, in general, the algorithm is faster and that it provides better results (cf. Fig.~\ref{Fig:k5OMP}). An expansion of the block algorithm to functions instead of vectors is straightforward. 
\section{Numerical Experiments}\label{Sec:Experiments}
The numerical experiments in the following will be restricted to representing two functions. The first target function will be the scattering solution of a plane wave on the surface on a sound hard cylinder. The second one will be the Hankel function of order 0 on the unit circle.
\subsection{Sound field on the cylinder}
\begin{figure}[!htb]
  \begin{center}
    \includegraphics[width=0.7\textwidth]{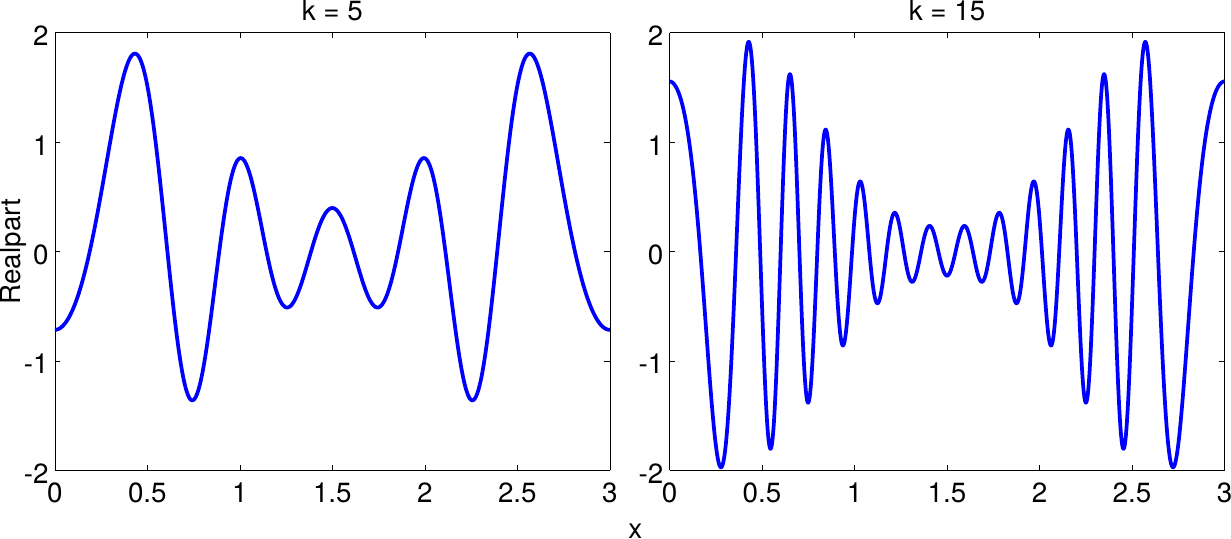}
    \caption{Real part of the solution on the cylinder cross section for $k = 5$ and $k = 15$.}\label{Fig:Cylinder}
  \end{center}
\end{figure}

As the first target function for the numerical experiments the description of the sound field on a sound hard circular cylinder caused a plane wave is used. This setup has the advantage that it can be reduced to a 2D problem and that an analytic solution is known (cf. \cite{CheWau07}):
\begin{equation}\label{Equ:CircleScat}
  f(\bx) = \frac{2}{\pi kr} \sum_{n=0}^\infty  \epsilon_n (-\text{i})^{(n-1)} \frac{\cos(n \phi)}{H'_n(kr)},
\end{equation}
where $\bx = re^{\I\phi}$, $\epsilon_0 = 1$, $\epsilon_n = 2, n>0$ and $H'_n$ is the derivative of the Hankel function of order $n$. The radius of the circle was set to $r = 1$. In Fig.~\ref{Fig:Cylinder} the real part of the target function $f(\bx)$ for wavenumbers $k = 5$ and $k = 15$ is shown, in both cases the sum in Eq.~(\ref{Equ:CircleScat}) was truncated after a fixed number of summands or when the norm of vector containing the values of the new summand on the whole circle was below a certain tolerance.

In the following the relative error between the analytic solution given in Eq.~(\ref{Equ:CircleScat}) and its expansion using a (for simplicity sampled) B-spline frame with window function $N_2(x)$ will be investigated. Of special interest will be the reconstruction error if only a small number of coefficients are used, e.g. 60 and 120 expansion coefficients. Since frames are used in the reconstruction different errors for different dual frames are expected. The interval of interest was set to $[0,3]$ and sampled using 601 points. The parameterization of the circle $\phi \in [0,2\pi]$ used in the analytic solution was re-scaled to this interval. This has the advantage that the used frame parameters $a = 1$ and $b = \frac13$ and the window function do not have to be re-scaled (see also Section \ref{Sec:RestrictionInterval}). 
The setup implies that 4 (shifted) windows are needed to cover the whole range of $[0,3]$. Although $f(x)$ is periodic, the B-spline window for generating the frame was not periodized. For calculating the frame coefficients using the dual frames, the interval of interest was extended to cover the full combined support of the necessary dual frames.
\begin{figure}[!h]
  \begin{center}
    \includegraphics[width=0.46\textwidth]{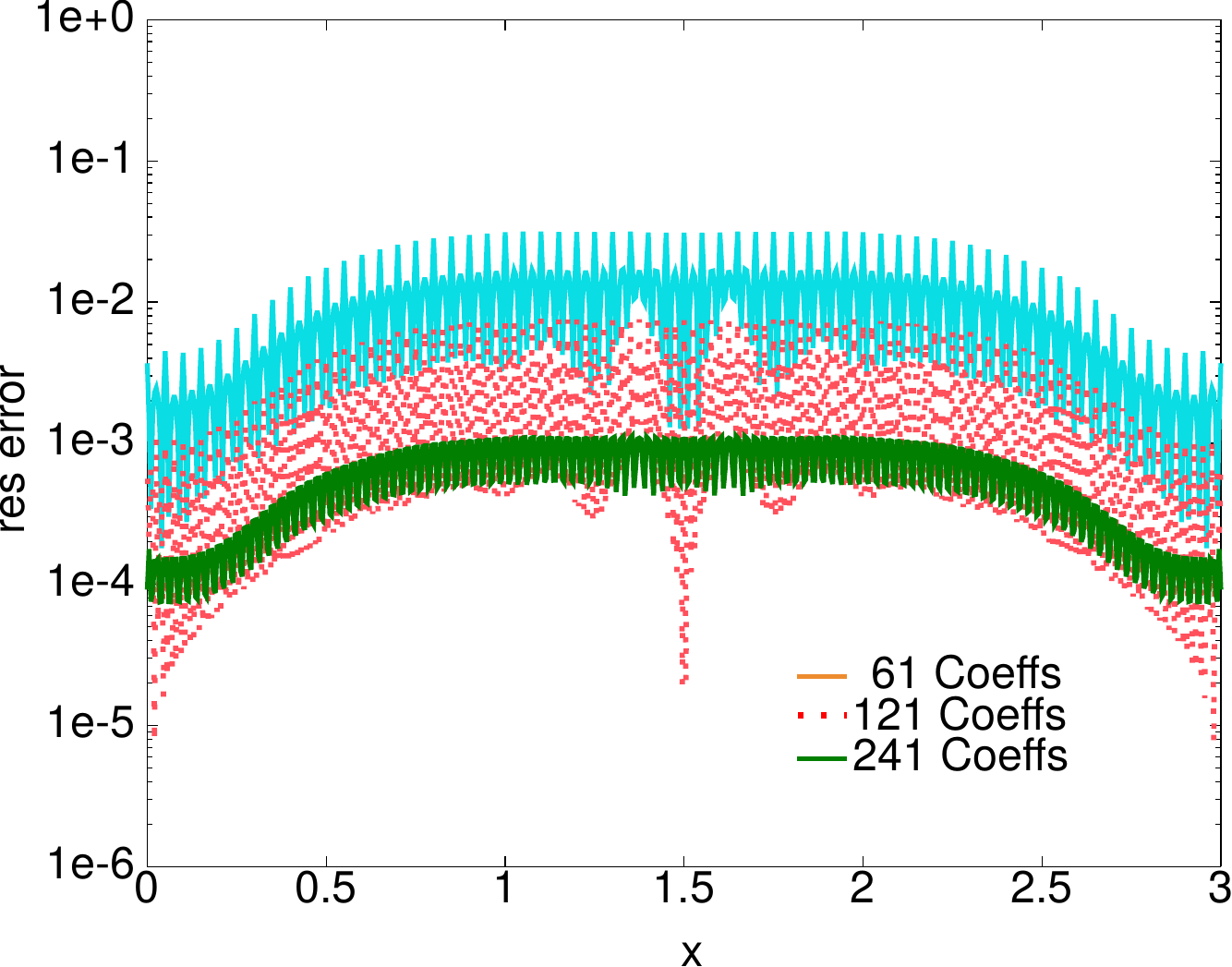}
    \includegraphics[width=0.46\textwidth]{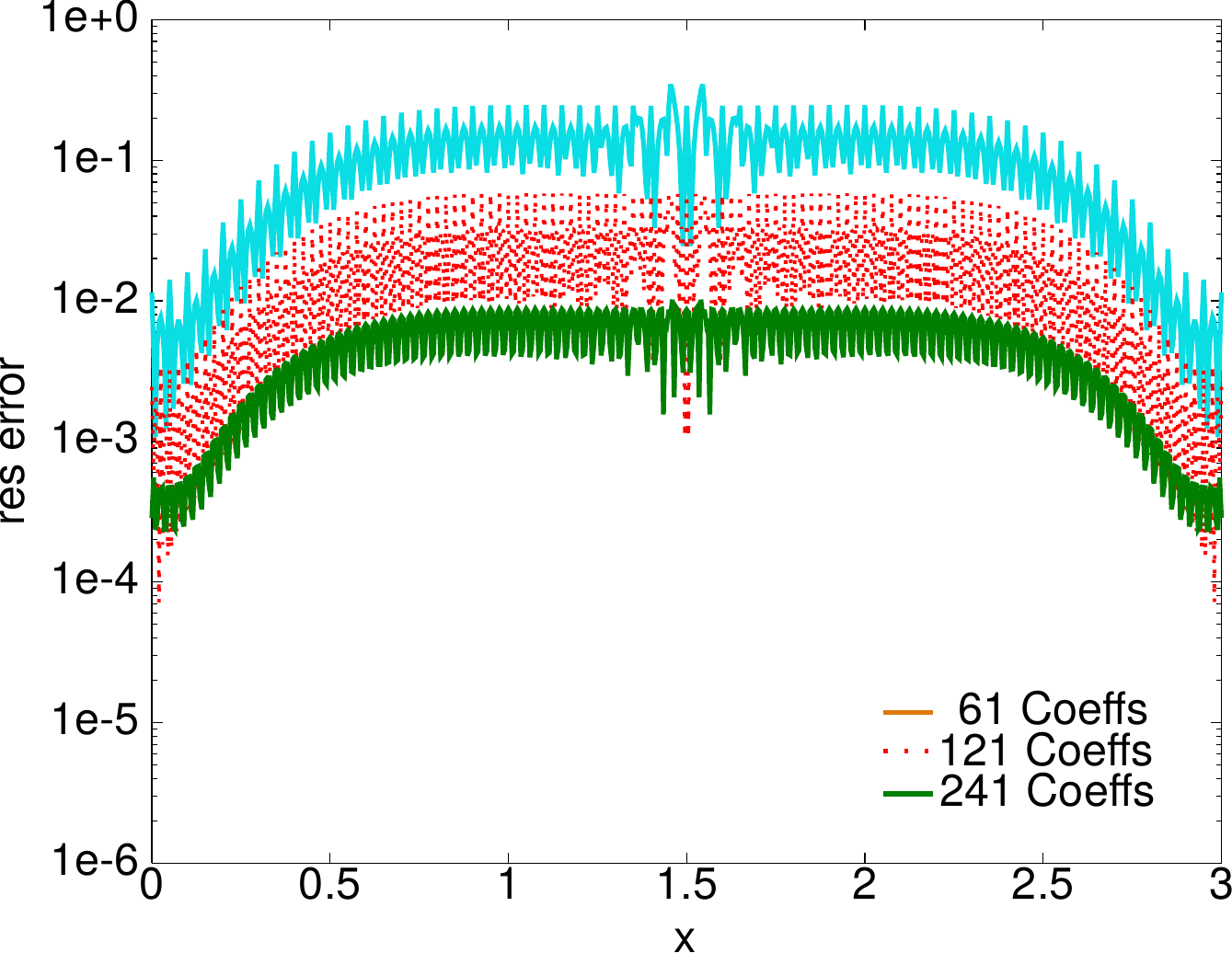}
  \caption{Relative error of the approximation of the target function using B-spline functions for the wavenumbers $k = 5$ (left side) and $k = 15$ (right side). The light blue/gray lines depict the approximation error when the function is expanded into a basis of 61 shifted versions of the scaled B-spline functions of order 2, the red dotted lines depict the approximation error when 121 shifted version of the B-spline are used, for the dark green lines 241 basis functions are used.}\label{Fig:ErrorBasis}
  \end{center}
\end{figure}

For reference the relative error between the original function for the wavenumbers $k = 5$ and $k = 15$ and the approximation using the standard B-spline ansatz functions based on the B-spline of order 2 with 61, 121, and 241 basis functions is given in Fig.~\ref{Fig:ErrorBasis}. This is equivalent to discretizing the circle with 60, 120, and 240 elements, and using standard (non-periodic and continuous) linear ansatz functions.  

For the numerical experiments with the (vector based) OMP method a slightly modified version of the MATLAB code from Stephen Becker \cite{Becker16} was used with three different blocksizes (1,10, and 20, see also Section \ref{Sec:OMP}).

\subsubsection{Wavenumber $k = 5$}\label{Sec:Wavek5}
\begin{figure}[!htb]
  \begin{center}
    \includegraphics[width=0.48\textwidth]{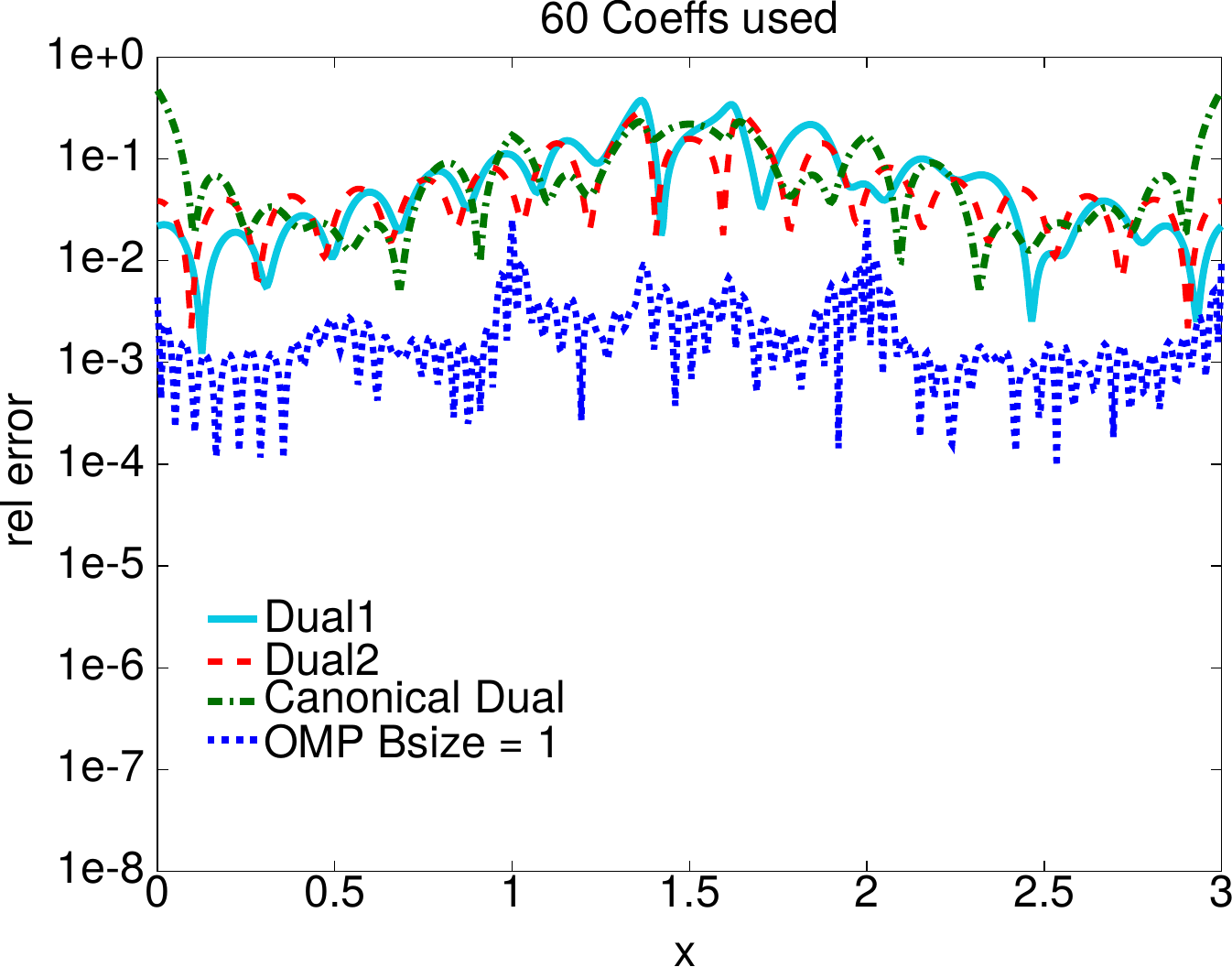}
    \includegraphics[width=0.48\textwidth]{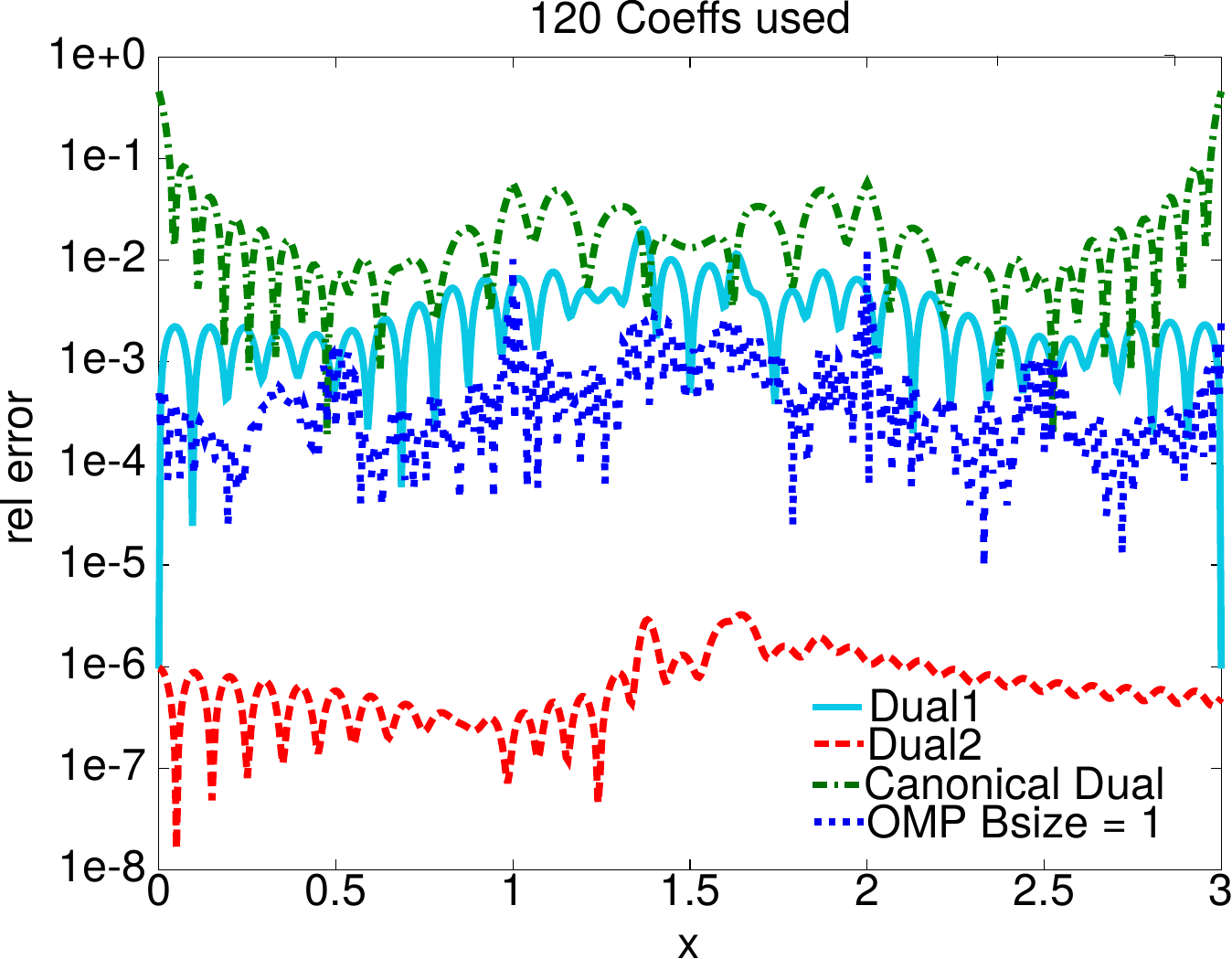}
    \caption{Left side: Relative approximation error when using the 60 biggest (in absolute value) coefficients using Dual1 (Eq.~(\ref{Equ:Dual1}), light blue line), Dual2 (red dashed line), the canonical dual (dark green dashed dotted line) and the OMP  algorithm with blocksize 1 (dotted dark blue line). Left side: Relative approximation error when using the 120 biggest (in absolute value) coefficients. In all cases $k = 5$.}\label{Fig:k5err60}
  \end{center}
\end{figure}

\begin{figure}[!htb]
  \begin{center}
    \includegraphics[width=0.48\textwidth]{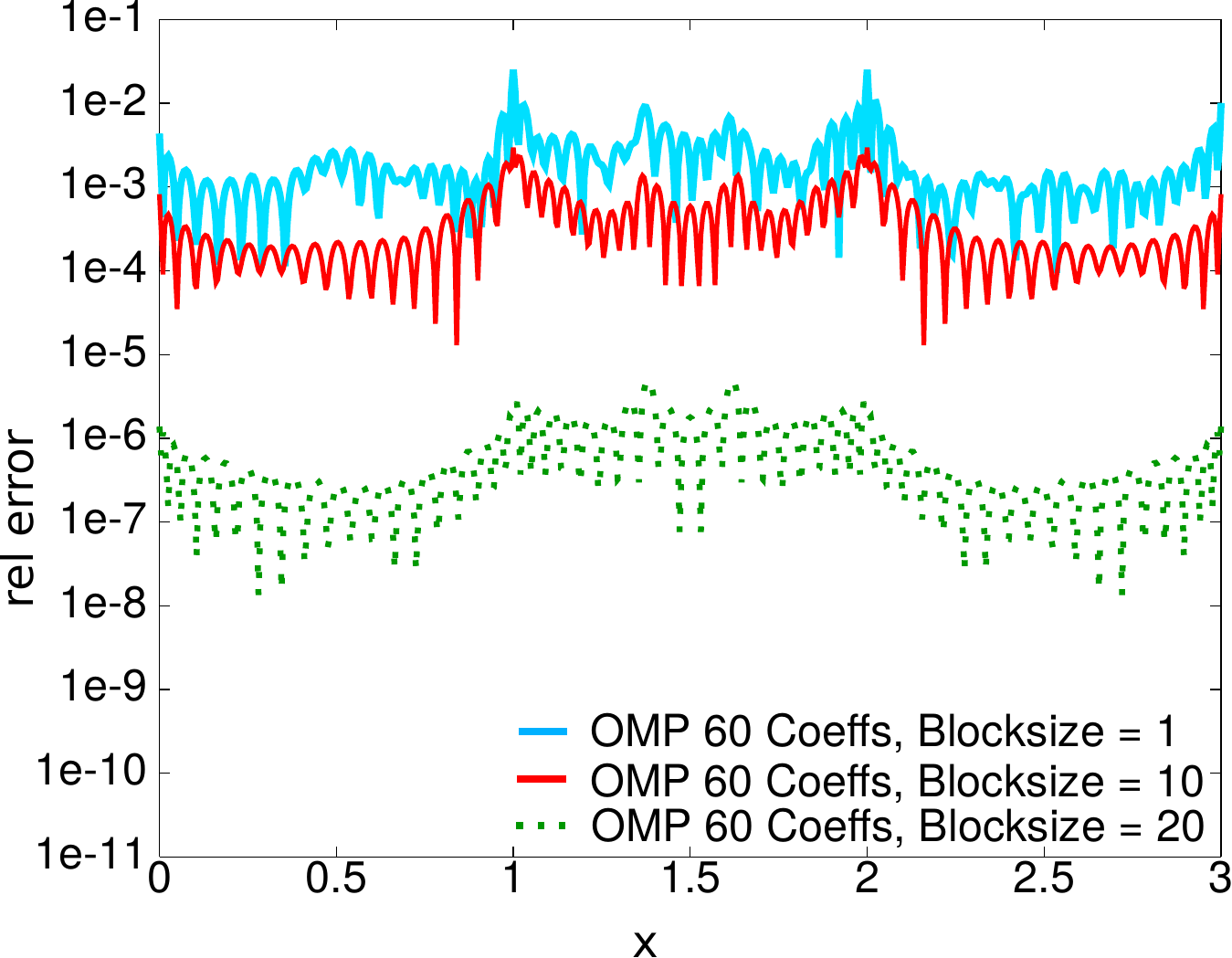}
    \includegraphics[width=0.48\textwidth]{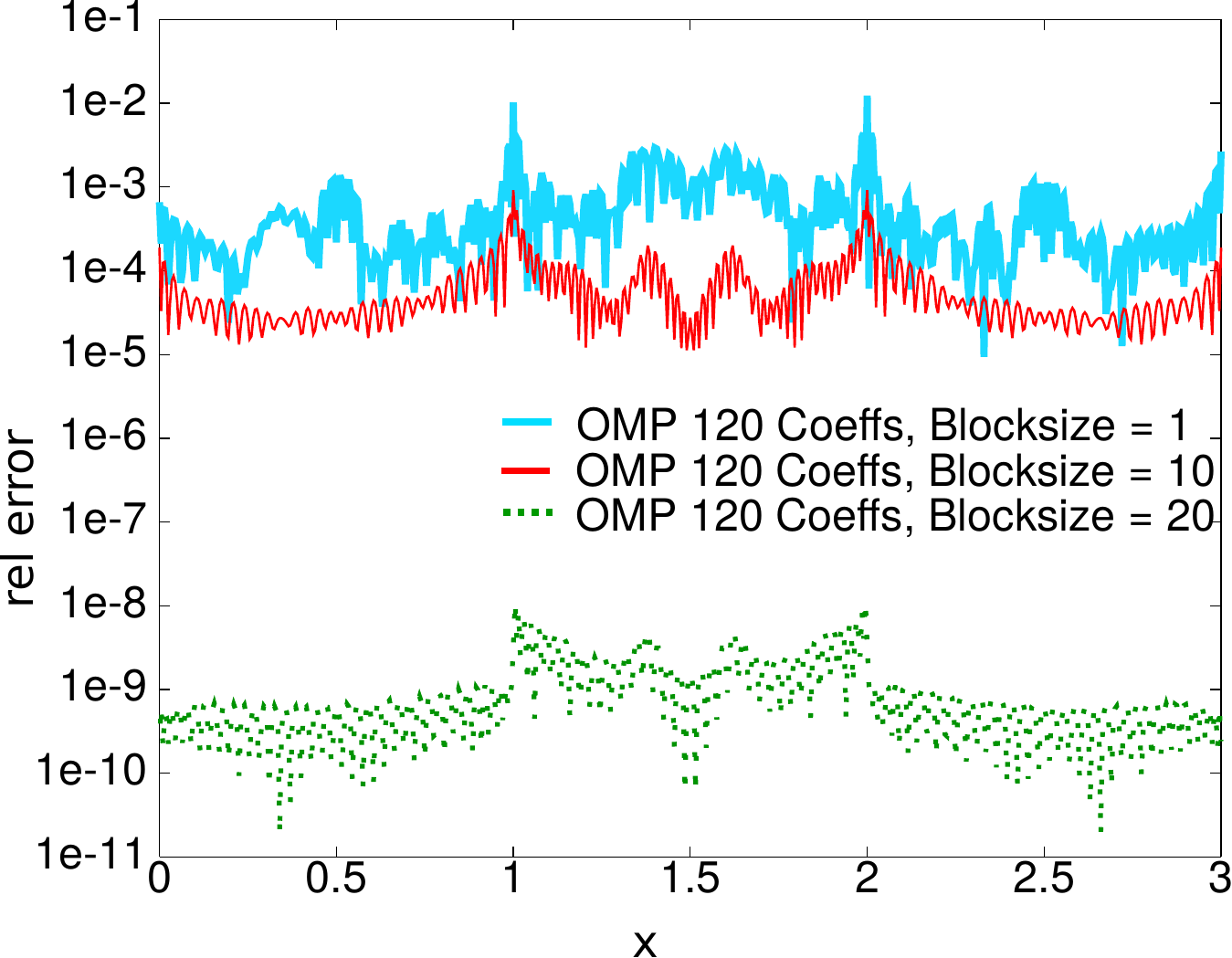}
    \caption{Relative error when the coefficients of the expansion are calculated using the OMP algorithm with different blocksizes. Left side: 60 coefficients are used. Right side: 120 coefficients are used. In all cases $k$ = 5.}\label{Fig:k5OMP}
  \end{center}
\end{figure}

In Fig.~\ref{Fig:k5err60} the errors of the representation of the total field on the cylinder/circle are given. In each figure the coefficients used in the expansion are calculated in 4 different ways. For the light blue lines the coefficients are calculated using the dual frame Dual1 defined in Eq.~(\ref{Equ:Dual1}), for the dashed red line the Dual2 described in Eq.~(\ref{Equ:Dual2}) is used, for the green dashed dotted line the canonical dual frame is used, and the dotted dark blue line is calculated using the OMP algorithm with blocksize $n = 1$. One the left side the largest 60 coefficients (in absolute value) were used for the representation, on the right side the largest 120 coefficients.%
Compared to the error when using the standard approach with a basis for piecewise linear functions (see Fig.~\ref{Fig:ErrorBasis}) the relative error when using the B-spline frame with only the 60 largest coefficients is relatively high. Only when using the OMP algorithm the errors are roughly in the same range. At a first glance this would point in the direction that at least for small wavenumbers and a small number of coefficients using more shifted versions of the original window has an advantage over modelling the oscillations using modulations of the window function. So does including oscillations in the ansatz functions of a frame pay off for low wavenumbers at all? After a closer look the answer is 'yes'. The redundancy of the frame offers more possibilities for approximating the target function. For example, when changing the search strategy in the OMP algorithm by increasing the blocksize the error can be reduced to a large degree as can be seen when looking at Fig.~\ref{Fig:k5OMP}, where the approximation error is depicted for blocksizes of 1, 10, and 20, respectively. If for example a blocksize of 20 is used the relative error can be reduced to a large degree. 
\begin{figure}[!htb]
  \begin{center}
    \includegraphics[width=0.8\textwidth]{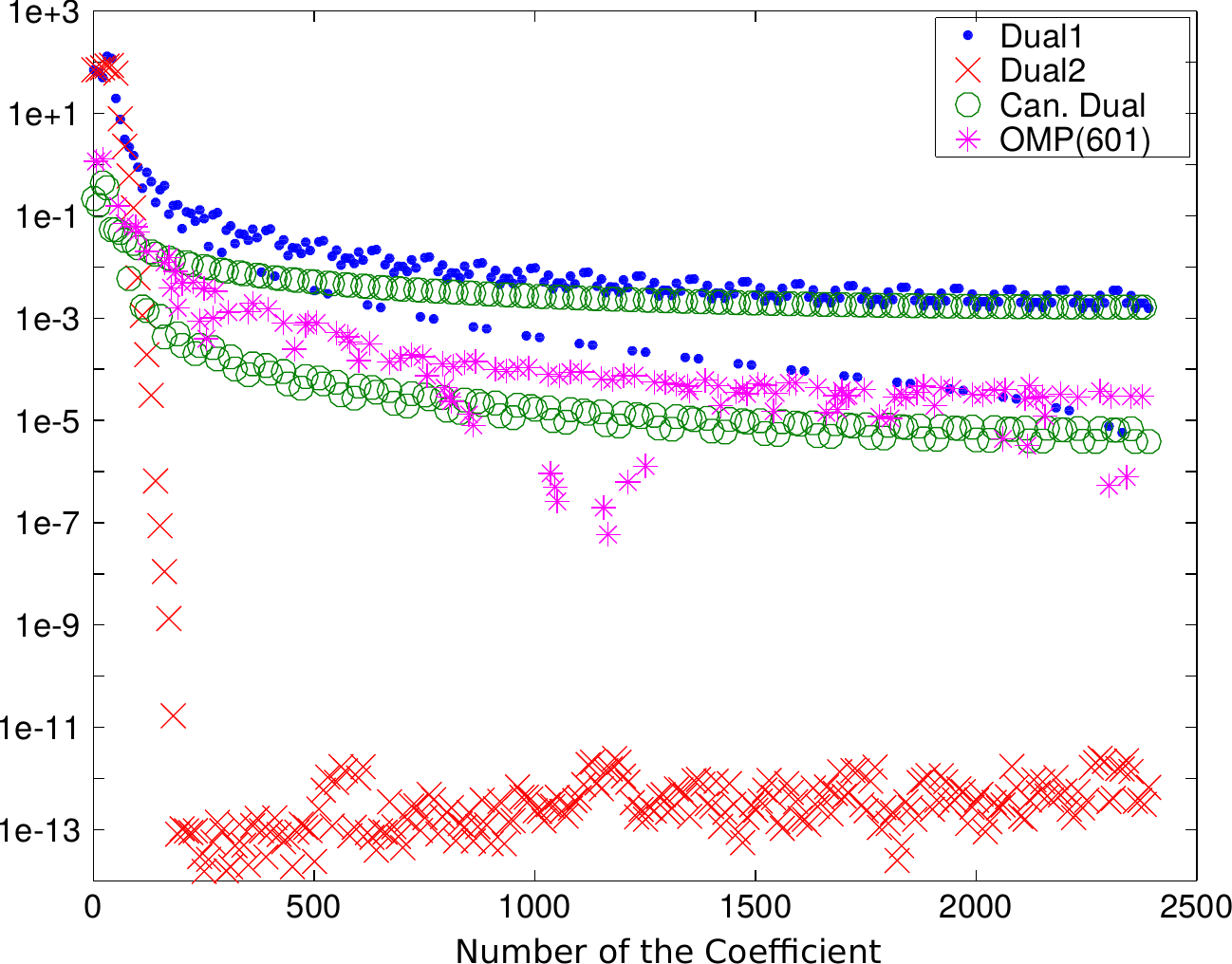}
    \caption{Absolute value of the coefficients used for the expansion of the target function. Coefficients marked with the dot were calculated with Dual1, coefficients marked with 'X' were calculated with Dual2, for the canonical dual and the OMP algorithm with blocksize = 1 'O' and '*' were used. Each dot in the figure represents one expansion coefficient, to keep the graph clearer only every 10th coefficient is plotted. The wavenumber $k$ = 5.}\label{Fig:k5Coeffs}
  \end{center}
\end{figure}

Redundancy in the frame means that there are several ways of approximating the target function, and each representation has its own properties. In Fig.~\ref{Fig:k5err60} the difference between different dual frames is clearly visible. Compared to the canonical dual and Dual1 the error for Dual2 decays very quickly with the number of coefficients used. This can also be seen in Fig.~\ref{Fig:k5Coeffs} that depicts the absolute value of the coefficients used for a complete reconstruction (relative errors ranging between $10^{-13}$ and $10^{-15}$ for the dual frames and $10^{-8}$ for the OMP with blocksize 1 using 601 coefficients and that was stopped when the absolute error was below a certain tolerance). The coefficients calculated with Dual2 decay very fast for the frame elements associated with a modulation factor bigger then $k = 5$, while the coefficients in all the other cases decay more slowly. For the canonical dual frame this fact is not so surprising because it can be associated with the pseudoinverse of the matrix containing the sampled frame elements, and it is known that using the pseudoinverse results in a solution that is, in general, not sparse but where the coefficient vector has the smallest $\ell_2$ norm. It is more surprising that the OMP algorithm with blocksize 1 that is designed to provide sparse solutions performs that badly compared to Dual2. 

But again, the performance of the OMP algorithm can be enhanced to a high degree if a blocksize bigger than one is used (see Sec.~\ref{Sec:OMP}). It can be seen in Fig.~\ref{Fig:k5OMP} that the performance of the OMP algorithm can be enhanced greatly if the blocksize is raised to 10 or 20, also the algorithm is much faster in these cases, because the cost for solving the small least square problems is small compared all the other calculations involved with the OMP.
\begin{figure}[!ht]
  \begin{center}
    \includegraphics[width=0.8\textwidth]{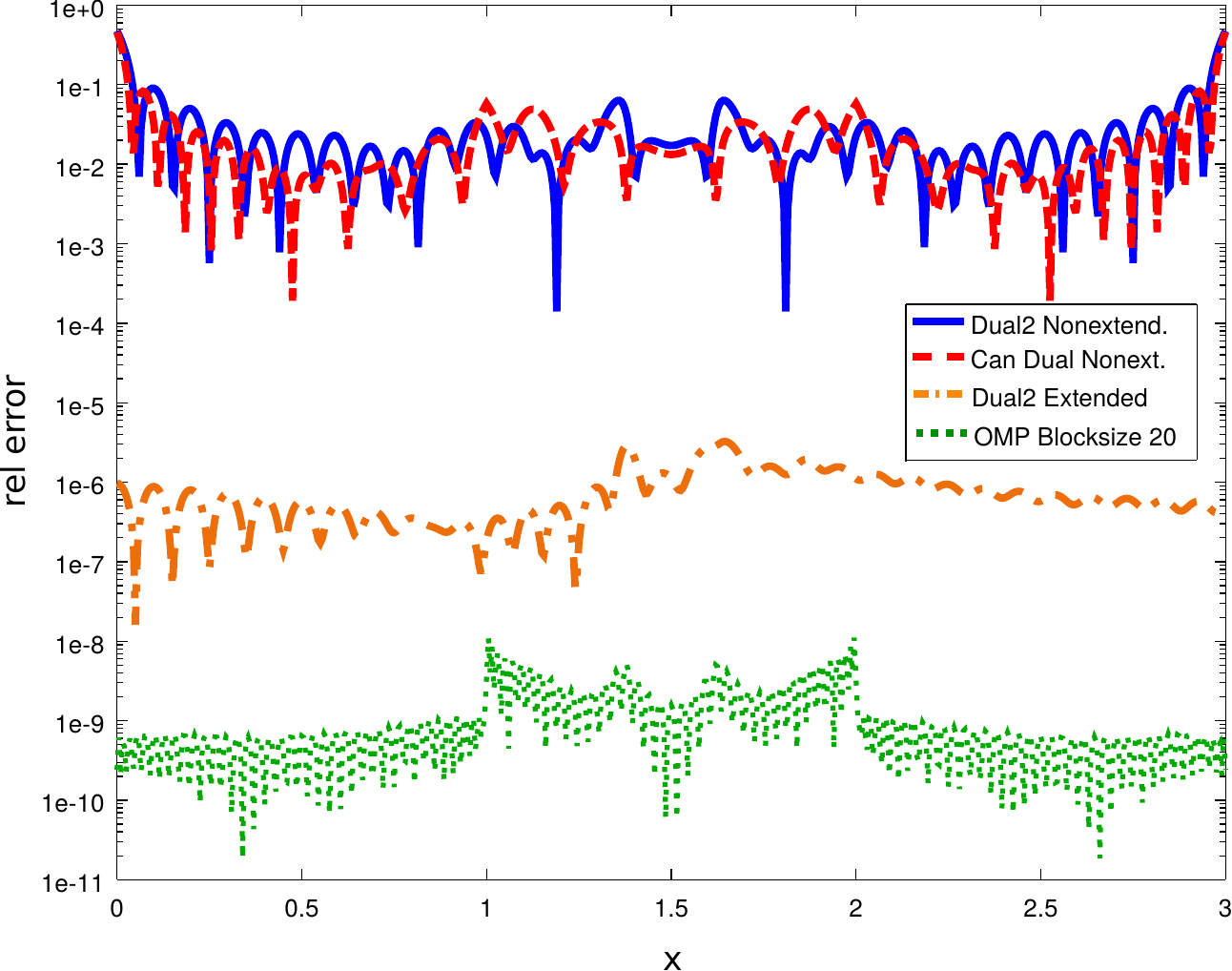}
    \caption{Relative error for the expansion with 120 frame coefficients calculated using Dual2 where the the interval of interest was \emph{not} extended (blue line), the canonical dual frame (dashed red line),  Dual2 where the interval was extended (dashed line) and the OMP algorithm with blocksize 20 (green dotted line).}\label{Fig:k5120CoeffNoExpand}
  \end{center}
\end{figure}

To illustrate the advantage of expanding the interval of interest Fig.~\ref{Fig:k5120CoeffNoExpand} shows the relative error for the expansion using 120 coefficients calculated with Dual2 \emph{without} expanding the interval of interest (blue line) compared with the canonical dual (dashed line), the Dual2 where the interval was extended (dashed dotted line) and the OMP algorithm with blocksize 20 (green dotted line). As the OMP does not rely on the dual, the error for the OMP is not dependent on the interval of interest.

\subsubsection{Wavenumber $k = 15$}

\begin{figure}[!htb]
  \begin{center}
     \includegraphics[width=0.8\textwidth]{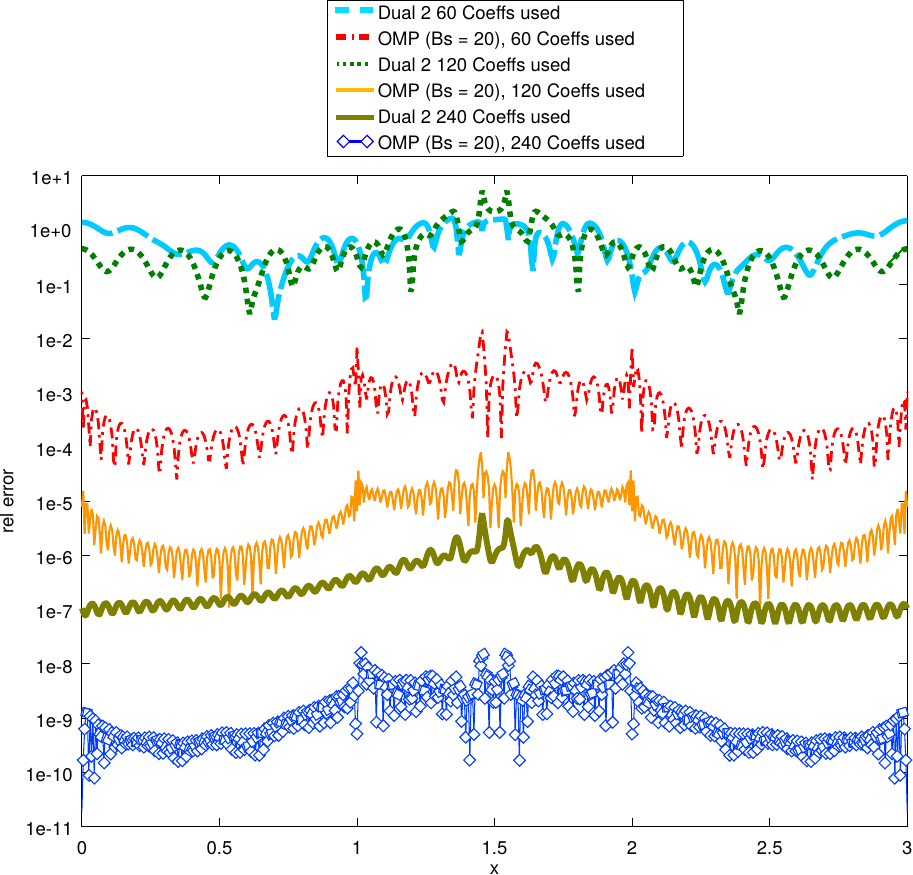}
  \caption{Relative approximation error for the wavenumber $k = 15$ using 60, 120, and 240 coefficients for Dual2 (light blue dashed line, dotted green line, army green thick line) and the OMP with blocksize 20 (dashed dotted red line, orange thin line, blue line with diamonds).}\label{Fig:k15err60}
  \end{center}
\end{figure}

For higher wavenumbers the block-OMP has a clear advantage when sparse solutions need to be found and the acceptable error tolerance is relatively high. Based on the experience gathered in Section \ref{Sec:Wavek5} the OMP algorithm was used with a blocksize of 20. When looking at the first column in Fig.~\ref{Fig:k15err60}, where only the 60 largest coefficients are used in the expansion, it can be observed that the errors for Dual2 is unacceptable high, the block OMP on the other hand provides an approximation where the maximum relative error is in the range of $10^{-4}$ to $10^{-2}$ which may be acceptable for some applications.
\begin{figure}[!htb]
  \begin{center}
    \includegraphics[width=0.7\textwidth]{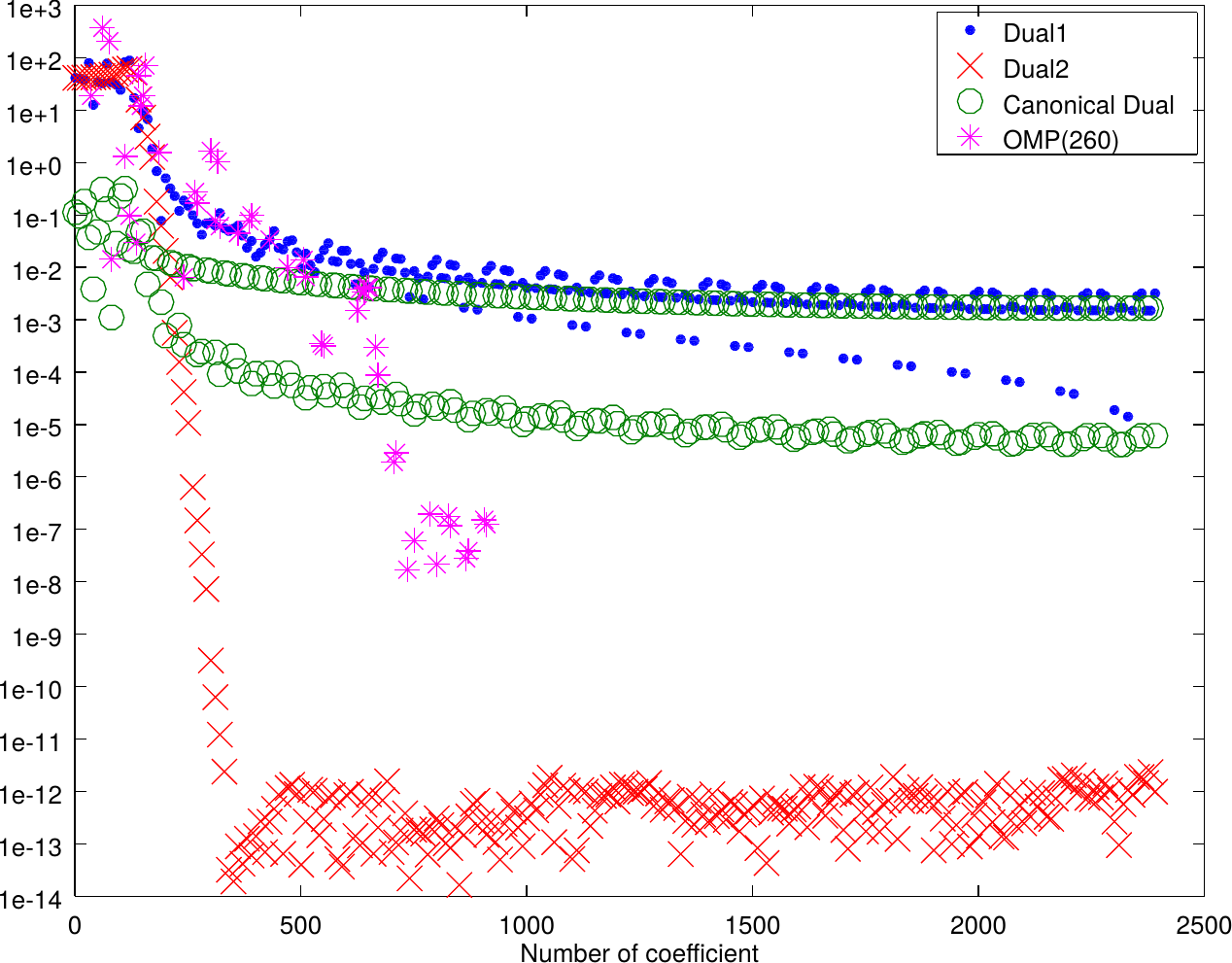}
    \caption{Absolute value of the coefficients for the representation of the target functions. Upper left: Dual1. Upper right: Coefficients are calculated using Dual2. Lower left: Canonical dual. Lower right: OMP with blocksize 20. In all cases the wavenumber was set to $k = 15$.}\label{Fig:k15coeffs}
  \end{center}
\end{figure}

In general, the behavior of the coefficients that are calculated using the dual frames is similar as for $k = 5$, however, when comparing Figs.~\ref{Fig:k5Coeffs} and \ref{Fig:k15coeffs} it becomes apparent that the coefficients calculated using Dual2 start to decay much later. For the case $k = 5$ only about 130 coefficients have absolute value above  $10^{-5}$, whereas for the case $k = 15$ about 240 are bigger than $10^{-5}$. This behavior is also reflected in Fig.~\ref{Fig:k15err60}. If only the 120 biggest coefficients (in absolute value, dotted green line ) calculated with Dual2 are used the approximation error is in the range of $10^{-1}$ to $10^1$, if the 240 largest coefficients are used the error (thicker army green line) is already below $10^{-5}$. In all cases the OMP algorithm with blocksize 20 again provides the best accuracy vs. sparsity relation.

\subsection{Hankel function $H_0$}
As second target function the Hankel function of order 0 was chosen. In 2D this function is the free space Green's function for the Helmholtz operator. $H_0(k||\bx - \bx^*||)$ describes the sound field at a point $\bx$ in 2D caused by a point source with wavenumber $k$ located at $\bx^*$. For the numerical experiment $\bx^* = (0,1.5)$, $k = 5$ and $k = 15$, and the field is evaluated on a circle with radius 1 and midpoint $(0,0)$.
\begin{figure}[!htb]
  \begin{center}
    \includegraphics[width=0.8\textwidth]{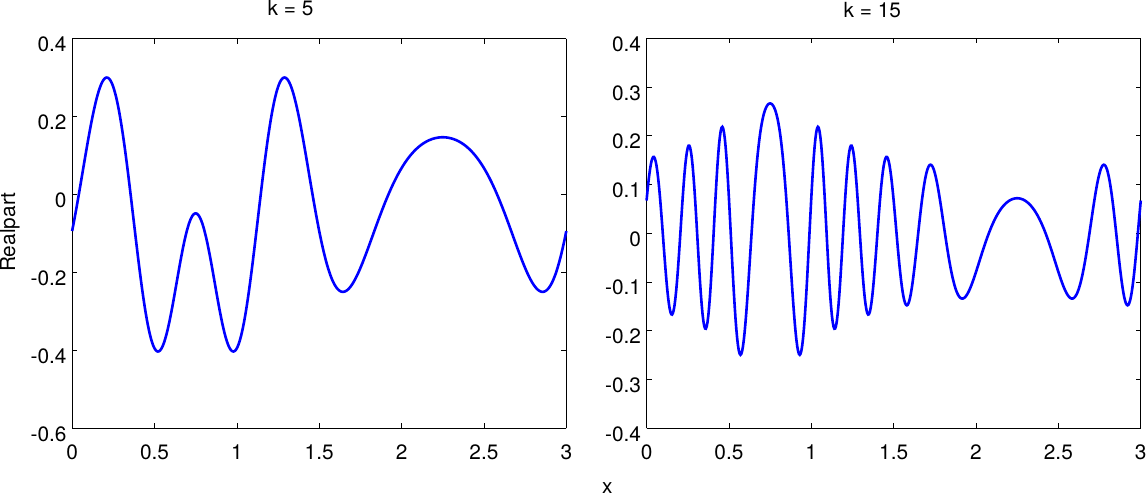}
    \caption{Real part of the acoustic field on a circle with radius 1 cause by a point source in 2D at (0,1.5) for wavenumbers $k = 5$ and $k = 15$.}\label{Fig:Hankel}
  \end{center}
\end{figure}
In Fig.~\ref{Fig:Hankel} the real part for this function is shown for wavenumbers $k = 5$ and $k = 15$.

For brevity we restrict ourselves to coefficients calculated using the dual frame Dual2 and the OMP algorithm with blocksize 20.
In general, the overall behaviour of the errors and the coefficients is similar to the case with the field on the cylinder. However, compared to the cylinder case, the decay of the coefficients calculated with Dual2 is slower. Naturally, this slower decay effects the size of the error, a comparison between the averaged errors is given in Tab.~\ref{Tab:Error}. For both cases using the block OMP results in smaller errors, however, it has to be added that first experiments also showed that the results for the OMP vary for different sound source position on a circle with radius $r = 1.5$, whereas in the Dual2 case, the errors always seem to be in the same range. 
\begin{table}[!h]
  \begin{tabular}{cccccccc}
    \hline
    \multicolumn{4}{c}{Target 1: Scattering solution} & \multicolumn{4}{c}{Target 2: Hankel function}\\
    \hline
     \# Coeffs & 60  & 120  & 240 & &  60  &   120  &  240 \\
    $k = 5$ &&&&&&&\\
    Dual2  & 6.3E-2 & 7.9E-7 & 3.6E-14 && 4.2E-2 & 3.2E-4 & 2.8E-7 \\
    OMP(20) & 7.6E-7 & 1.2E-9 & 3.5E-12 && 2.5E-4 & 3.0E-11 & 1.0E-13\\
    \hline
    $k = 15$ &&&&&&&\\
    Dual2  & 6.3E-1 & 5.5E-1 & 4.0E-7 && 4.7E-1 & 9.7E-2 & 5.7E-6\\
    OMP(20) & 9.0E-4 & 7.0E-6 & 1.9E-9 && 3.5E-3 & 2.8E-5 & 3.4E-13\\
  \hline
  \end{tabular}
  \caption{Comparison of the average relative errors for approximating the target functions using coefficients calculated with Dual2 and the OMP algorithm with blocksize 20. Both functions are approximated for wavenumbers $k = 5$ and $k = 15$ using 60, 120 and 240 coefficients.}\label{Tab:Error}
\end{table}

\section*{Summary and Outlook}
\begin{figure}[!htb]
  \begin{center}
    \includegraphics[width=0.7\textwidth]{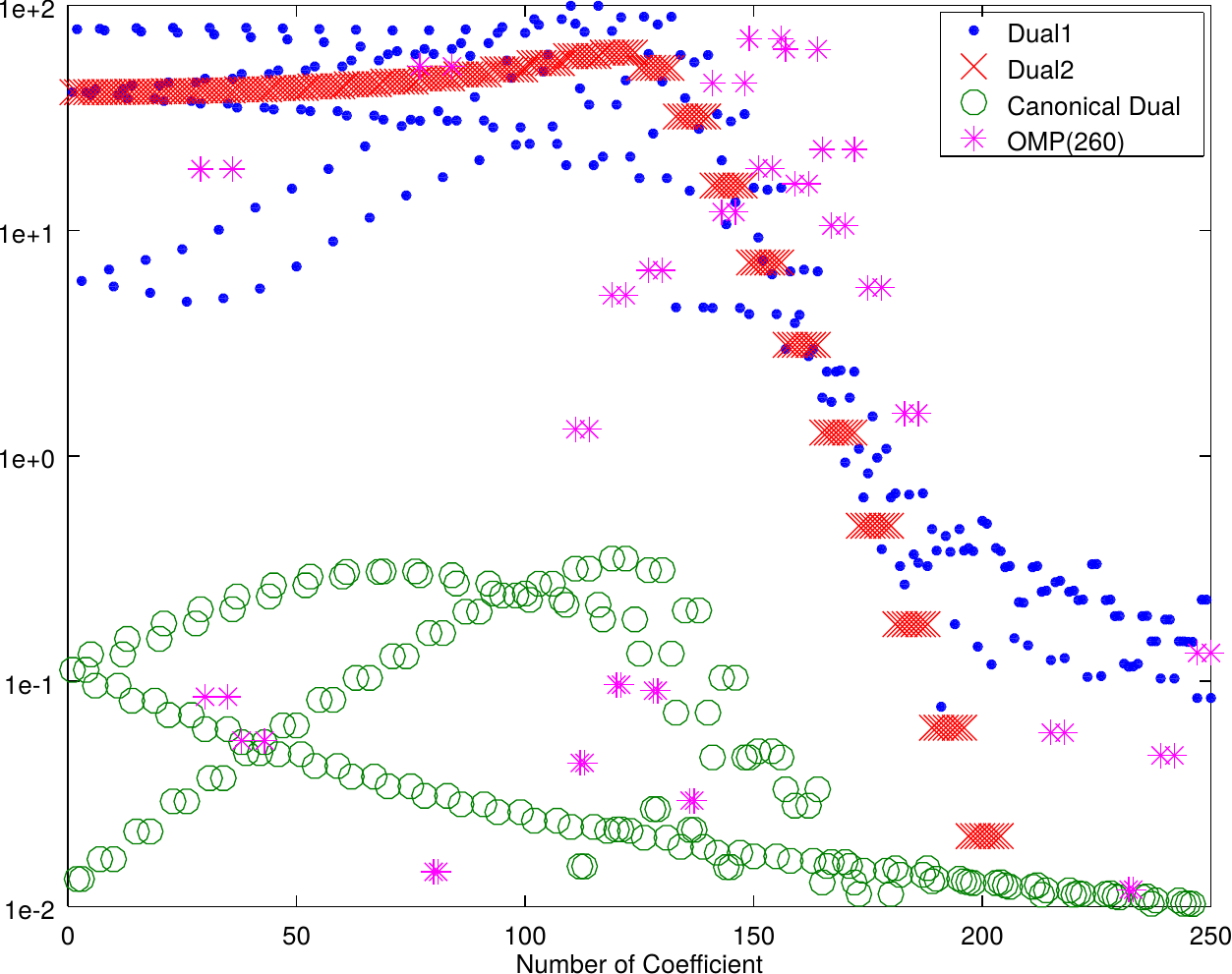}
    \caption{Zoomed version of Fig.~\ref{Fig:k15coeffs}.}\label{Fig:Coeffszoomed}
  \end{center}
\end{figure}

In conclusion, numerical experiments show that Gabor frames based on B-splines can be used to efficiently represent solutions of the Helmholtz equation, especially if methods used in compressed sensing are used to calculate the coefficients. If B-spline functions are used as generating functions for the Gabor frame, easy ways can be provided to find frame parameters and to construct dual frames, however, there are several ways for finding the expansion coefficients that have different properties. The redundancy of the frame offers several possibilities to construct approximations, and a good way of finding the right frame coefficients needs to be found. Using methods like the OMP algorithm have certain advantages if a good strategy of finding the right candidates (see Step 2 in Section~\ref{Sec:OMP}) can be developed and implemented. When looking for example at Fig.~\ref{Fig:Coeffszoomed} that contains the zoomed version of Fig.~\ref{Fig:k15coeffs} it becomes apparent that the largest coefficients calculated using Dual2 are concentrated around the frame elements with modulations related to the wavenumber. Thus besides using blocksizes bigger than 1, one could in the first few steps concentrate the search for the right candidates to the frame elements with modulations close to the wavenumber of interest, and extend then the search interval to a bigger set of frame elements.

For future work it will also be necessary to look at problems in 3D. On the one hand the extension to 3D is straight forward by using tensor products of one dimensional frames, on the other hand it will be important and necessary to adapt the frames to provide stiffness matrices that can be easily approximated by sparse matrices or that can be used in connection with methods like the fast multipole method or $\mathcal{H}$-matrices. To that end it will be necessary to adapt the frames to the Green's function for the Helmholtz equation in 3D, and to provide efficient quadrature methods. 

\section*{Acknowledgements}
This study was supported by the Austrian Science Fund (FWF), Project Biotop [I 1018-N25].

\appendix
\section{OMP for functions}\label{Sec:OMPL2}
To apply the formulation of the OMP algorithm used in Section \ref{Sec:OMP} to complex valued functions some modifications need to be made, which results in a sort of weak formulation. Since the residuum cannot be addressed directly, one needs to look at ${\bf r}'_i = \langle r_n(x), g_i(x) \rangle$, where $r_n(x) = f(x) - \sum_{j=1}^n \gamma_j g_j(x)$ is the residuum after the $n$-th step. 
\begin{enumerate}
\item ${\bf r}_i' = \langle f, g_i \rangle \; \text{for } i = 1,\dots N$
\item Find the right candidate $i_0 = \text{argmax} |{\bf  r'}|$
\item Add the candidate to the list of used frame atoms: $\mathcal{I} = \mathcal{I} \cup \{i_0\}$
\item Find ${\bf \gamma} = \text{argmin}_{\bf c} || f - \sum\limits_{i\in\mathcal{I}} c_i g_i ||^2$, where $||.||$ denotes the norm in $L^2(\mathbb{R})$.
\item Update ${\bf r}_i' = \left\langle  f - \sum\limits_{j\in\mathcal{I}}\gamma_j g_j, g_i \right\rangle$
  
\item If $\text{number}_\text{iterations} < \text{max}_\text{iterations}$ jump to step 2,
\end{enumerate}
Again, instead if using the number of iterations as stopping criterion, it is also possible to stop the iterations if $||{\bf r}'||_2$ (now in $\mathbb{C}^N)$ is below a tolerance provided by the user.

To find the minimum in Step 4 
\begin{equation}\label{Equ:Lproblem}
  \min_ {\bf c}  ||f(x) - \sum\limits_{i\in\mathcal{I}} c_i g_i(x)||^2 = \min_{\bf c}  \left\langle f(x) - \sum\limits_{i\in\mathcal{I}} c_i g_i(x), f(x) - \sum\limits_{i\in\mathcal{I}} c_i g_i(x) \right\rangle
\end{equation}
the vector product is calculated and the derivatives with respect to the real and imaginary part of each $c_i$ need to be set to 0.

The product in Eq.~\ref{Equ:Lproblem} can be expanded into
\begin{equation}\label{Equ:ComplxMult}
  \langle f,f \rangle -
  \sum_{i\in\mathcal{I}} \langle f,g_i \rangle c^*_i - \sum_{i\in\mathcal{I}} \langle g_i,f \rangle c_i +
  \sum_{i,j \in \mathcal{I}} c_i \langle g_i, g_j \rangle {c^*_j},
\end{equation}
where $c_i^*$ the complex conjugate of the complex number $c_i$.

Taking the derivatives with respect to the real and imaginary parts of each $c_i, i\in\mathcal{I}$ results in the two following equations:
\begin{align}
   \sum_{j\in\mathcal{I}}\langle g_i,g_j \rangle c^*_j + \sum_{j\in\mathcal{I}} \langle g_j,g_i\rangle c_j - \langle f,g_i \rangle - \langle g_i,f \rangle &= 0,\\
   \sum_{j\in\mathcal{I}} \langle g_i,g_j \rangle c^*_j - \sum_{j\in\mathcal{I}} \langle g_j,g_i\rangle c_j + \langle f,g_i \rangle - \langle g_i,f \rangle  &= 0,
\end{align}
which yields the linear system ${\bf G^T c} = {\bf f}$ with ${\bf G^T}$ being the transpose of the matrix ${\bf G}$ defined by ${\bf G}_{ij} = \langle g_i,g_j \rangle$, ${\bf f}_i = \langle f,g_i \rangle$, and ${\bf c}$ contains the complex valued vector of the unknown coefficients. Note, that the sizes of ${\bf G}$ and ${\bf f}$ only depend on the size of the (small) set $\mathcal{I}$.
\section*{References}
\bibliography{kreuzer18}
\end{document}